\documentclass{article}
\usepackage{amsmath,amssymb,mathrsfs,amsthm,verbatim, MnSymbol, fancyhdr}

\newtheorem{Theorem5}{Theorem}

\newtheorem{Lemma1}{Lemma}[section]
\newtheorem{Lemma2}[Lemma1]{Lemma}
\newtheorem{Lemma3}[Lemma1]{Lemma}
\newtheorem{Lemma4}[Lemma1]{Lemma}
\newtheorem{Lemma5}[Lemma1]{Lemma}
\newtheorem{Lemma6}[Lemma1]{Lemma}
\newtheorem{Corollarymono}[Lemma1]{Corollary}

\date{}

\begin{document}

\title{Seiberg-Witten Flow in Higher Dimensions}
 \author{ {Lorenz Schabrun} \footnote{
Department of Mathematics, The University of Queensland, Brisbane,
QLD 4072, Australia. Email: lorenz@maths.uq.edu.au}}

\maketitle \textbf{Abstract.} We show that for manifolds of dimension $m\geq5$, the flow of a Seiberg-Witten-type functional admits a global smooth solution.

\section{Introduction}
The Seiberg-Witten invariant has proven a very effective tool in four-dimensional geometry. Its computation involves finding non-trivial solutions to the system of first order Seiberg-Witten equations, called monopoles. Monopoles represent the zeros of the Seiberg-Witten functional (\ref{swfunc}) (see \cite{Jost.Variational}). In \cite{Schabrun.Flow}, the flow for the Seiberg-Witten functional on a 4-manifold was studied. It was shown that the flow admits a global solution which converges in $C^\infty$ to a critical point of the functional.

The Seiberg-Witten equations and functional do not generalize immediately to higher dimensions, since they depend on the notion of self-duality on four dimensional manifolds. Nonetheless, a number of generalizations of Seiberg-Witten theory have been suggested for higher dimensional manifolds, see for example \cite{SWHigher.Tian}, \cite{SWHigher.Bilge} and \cite{SWHigher.G2}. In this paper, we extend the global existence result obtained for the Seiberg-Witten functional in \cite{Schabrun.Flow} for dimension 4 to a functional of similar form in higher dimensions.

The main additional estimate in establishing global existence in higher dimensions is a so-called monotonicity formula. This idea was used by Struwe for the heat flow of harmonic maps in higher dimensions \cite{St1}, and has also been used to study the Yang-Mills and Yang-Mills-Higgs flows in higher dimensions, see \cite{CS} and \cite{Hong.Tian}. See also \cite{Lingradient} and \cite{LinWang} for the harmonic map flow, and \cite{TianCalibrated} for sequences of weakly converging Yang-Mills connections.

Let $M$ be a compact oriented Riemannian $m$-manifold which admits a
$\operatorname{Spin^c}$ structure $\mathfrak s$. Denote by $\mathcal S= W\otimes
\mathcal L$  the corresponding spinor bundle and by $\mathcal S^ \pm=
W^{\pm}\otimes \mathcal L$ the half spinor bundles, and by $\mathcal
L^2$ the corresponding determinant line bundle. Recall that the
bundle $\mathcal S^+$ has fibre $\mathbb{C}^2$. Let $A$ be a
unitary connection on $\mathcal L^2$. Note that we can write
$A=A_0 + a$, where $A_0$ is some fixed connection and $a \in
i\Lambda ^1 M$ with $i=\sqrt {-1}$. Denote by $F_A = dA \in
i\Lambda ^2 M$ the curvature of the line bundle connection $A$.
Let $\{e_j\}$ be an orthonormal basis of $\mathbb R^4$.  A
$\operatorname{Spin(4)^c}$-connection on the bundles $\mathcal S$ and $\mathcal
S^ \pm$ is locally defined by
\begin{equation} \label{connectiondef}
\nabla _A  = d + \frac{1}{2}(\omega +A),
\end{equation}
where  $\omega=\frac 1 2 \omega_{jk}e_je_k$ is induced by the
Levi-Civita connection matrix  $\omega_{jk}$ and $e_j e_k$ acts by
Clifford multiplication (see \cite{Jost.Riemannian}). We denote the curvature of $\nabla _A$ by $\Omega _A$. 
We define the configuration space $\;\Gamma (\mathcal S^ +  )
\times \mathscr A$, where $\mathscr A$ is the space of unitary
connections on $\mathcal L^2$, and let $(\varphi ,A) \in \Gamma
(\mathcal S^ +  ) \times \mathscr A$.

We first recall the definition of the Seiberg-Witten functional on 4-manifolds. The Seiberg-Witten functional $\mathcal{SW}:\;\Gamma (\mathcal S^
+  ) \times \mathscr A \to \mathbb{R}$ is given by
\begin{equation} \label{swfunc}
\mathcal{SW}(\varphi ,A) = \int_M {\left| {\nabla _A \varphi }
\right|^2  + \left| {F_A^ +  } \right|^2  + \frac{S} {4}\left|
\varphi  \right|^2  + \frac{1} {8}\left| \varphi  \right|^4 }\,dV,
\end{equation}
where $S$ is the scalar curvature of $M$. The Seiberg-Witten functional (\ref{swfunc}) is invariant under the action of a gauge group. The group of gauge transformations is
\[
\mathscr G = \left\{ {g:M \to U(1)} \right\}.
\]
$\mathscr G$ acts on elements of the configuration space via
\[
g^*(\varphi ,A) = (g^{-1}\varphi ,A + 2g^{-1}dg).
\]
Using the relation
\begin{equation} \label{selfdualrelation}
\left\| {F_A } \right\|_{L^2 }  = 2\left\| {F_A^ +  } \right\|_{L^2 }  - 4\pi ^2 c_1 (\mathcal L^2)^2, 
\end{equation}
where $c_1 (\mathcal L)$ is the first Chern class of $\mathcal{L}$ (see \cite{wildworld}), one can also write the functional in the form
\begin{equation} \label{swfuncchern}
\mathcal{SW}(\varphi ,A) = \int_M {\left| {\nabla _A \varphi }
\right|^2  + \frac{1}{2}\left| {F_A} \right|^2  + \frac{S} {4}\left|
\varphi  \right|^2  + \frac{1} {8}\left| \varphi  \right|^4 }\,dV + \pi ^2 c_1 (\mathcal L)^2.
\end{equation}
Now, consider again the case of an $m$-manifold $M$. The functional (\ref{swfunc}) is not defined here, since self-duality is a phenomenon that occurs only in dimension 4. However, we may use (\ref{swfuncchern}) to extend the Seiberg-Witten functional to higher dimensions. Note that the constant term $\pi ^2 c_1 (\mathcal L)^2$ does not affect the Euler-Lagrange equations and so is irrelevant for the results in this paper. The Euler-Lagrange equations for the Seiberg-Witten functional are
\begin{equation} \label{eq1}
 - \nabla _A^* \nabla
_A \varphi  - \frac{1}{4}\left[ {S + \left| \varphi  \right|^2 }
\right]\varphi =0,
\end{equation}
\begin{equation} \label{eq2}
 - d^* F_A   - i\operatorname{Im} \left\langle {\nabla_A \varphi ,\varphi }
\right\rangle =0,
\end{equation}
and as in \cite{Schabrun.Flow}, we define the flow of the Seiberg-Witten functional by
\begin{equation} \label{flow1}
\frac{{\partial \varphi }} {{\partial t}} =  - \nabla _A^* \nabla
_A \varphi  - \frac{1}{4}\left[ {S + \left| \varphi  \right|^2 }
\right]\varphi,
\end{equation}
\begin{equation} \label{flow2}
\frac{{\partial A}} {{\partial t}} =  - d^* F_A   -
i\operatorname{Im} \left\langle {\nabla_A \varphi ,\varphi }
\right\rangle
\end{equation}
with initial data
\[
(\varphi (0),A(0)) = (\varphi _0 ,A_0).
\]
Regarding the existence of solutions to the flow (\ref{flow1}) and (\ref{flow2}), we prove the following theorem.

\begin{Theorem5} \label{globalexisthigher}
For any given smooth  $(\varphi_0, A_0)$ and $m$-dimensional Riemannian manifold $M$ for $m \geq 5$, equations (\ref{flow1})
and (\ref{flow2}) admit a unique global smooth solution on $M
\times \left[ {0,\infty } \right)$ with initial data $(\varphi_0,
A_0)$.
\end{Theorem5}

As in \cite{Schabrun.Flow}, we have the energy inequality
\begin{equation} \label{energyinequality1}
\frac{d} {{dt}}\mathcal{SW}(\varphi (t),A(t)) =  - \int_M {\left[ 2{\left|
{\frac{{\partial \varphi }} {{\partial t}}} \right|^2  + \left|
{\frac{{\partial A}} {{\partial t}}} \right|^2 } \right]}
\leqslant 0.
\end{equation}
or
\begin{equation} \label{energyinequality2}
\int_0^T  {\left[ 2{\left\| {\frac{{\partial \varphi }} {{\partial
t}}} \right\|_{L^2 }^2  + \left\| {\frac{{\partial A}} {{\partial
t}}} \right\|_{L^2 }^2 } \right]}  =  \mathcal{SW}(\varphi_0,A_0)
- \mathcal{SW}( \varphi (T),A(T) ).
\end{equation}
Many of the proofs in \cite{Schabrun.Flow} do not contain dimensional considerations, and are also valid in the $m$-dimensional case. For the proofs of the following Lemmas, we direct the reader to that paper.

\begin{Lemma1} \label{Local}
For any given smooth initial data $(\varphi_0, A_0)$, the system (\ref{flow1})-(\ref{flow2}) admits a unique local smooth
solution on $M \times \left[ 0,T \right)$ for some $T>0$.
\end{Lemma1}

\begin{Lemma5} \label{phibound}
Let $(\varphi,A)$ be a solution of (\ref{flow1})-(\ref{flow2})
on $M \times [0,T)$, and write $m=\mathop {\sup }\limits_{x \in M}
\left| {\varphi _0 } \right|$. Then for all $t\in [0,T)$, we have
\begin{equation} \label{varphibound}
\mathop {\sup_{x\in M} } \left| {\varphi} (x,t) \right| \leqslant
\max \{ m,\sqrt {\left| {S_0 } \right|} \}.
\end{equation}
\end{Lemma5}

\begin{Lemma2} \label{firstderivativeestimate}
There exist positive constants $c,c'$ such that the following
estimate holds:
\begin{align*}
&\frac{\partial } {{\partial t}}\left( {\left| {\nabla _A \varphi
} \right|^2  + \left| {F_A } \right|^2 } \right) + \Delta \left(
{\left| {\nabla _A \varphi } \right|^2  +\left| {F_A } \right|^2 }
\right)  \\
&\leqslant - c'\left( {\left| {\nabla _A^2 \varphi } \right|^2  +
\left| {\nabla F_A } \right|^2 } \right) + c\left( {\left| {F_A }
\right| + 1} \right)\left( {\left| {\nabla _A \varphi } \right|^2
+ \left| {F_A } \right|^2 +1} \right).
\end{align*}
\end{Lemma2}

\begin{Lemma6} \label{higherderivativeestimate} Let $(\varphi, A)$ be a solution to (\ref{flow1})-(\ref{flow2})  in $M\times [0,T)$
with initial values $(\varphi_0 ,A_0)$. Suppose $\left| {\nabla _A
\varphi } \right| \leqslant K_1$ and $\left| F_A   \right|
\leqslant K_1$ in $M\times [0, T)$ for   some constant $K_1>0$.
Then for any positive integer $n\geq 1$, there is a constant
$K_{n+1}$ independent of $T$  such that
\[\left|
{\nabla _A^{(n+1)} \varphi } \right| \leqslant K_{n+1},\quad
\left| {\nabla _M^{(n)} F_A } \right| \leqslant K_{n+1}\quad \mbox
{in } M\times [0, T),
\]
where $(n)$ denotes $n$ iterations of the derivative.
\end{Lemma6}

\section{Estimates}

We first derive a monotonicity inequality for the flow (\ref{flow1})-(\ref{flow2}). We define
\[
e(\varphi,A)(x,t)=\left| \nabla_A \varphi \right|^2+\frac{1}{2} \left| F_A \right|^2 + \frac{S}{4} \left| \varphi \right|^2 + \frac{1}{8} \left| \varphi \right|^4.
\]
Let $z=(x,t)$ denote a point of $M \times \mathbb{R}$, with $z_0=(x_0,t_0) \in M \times [0,T]$. We define
\[
T_R (z_0 ) = M \times \left[ {t_0  - 4R^2, t_0  - R^2 } \right],
\]
and
\[
P_R (z_0 ) = B_R (x_0 ) \times \left[ {t_0  - R^2 ,t_0 } \right],
\]
where $B_R(x_0) \subset M$ denotes a ball centered at $x_0$ with radius $R$. Note that in constucting $T_R(z_0)$ we require that $t_0 -4R^2 \geq 0$ or $R \leq \sqrt{t_0} /2$. We abbreviate $T_R(0,0)=T_R$ and $P_R(0,0)=P_R$. The fundamental solution to the backward heat equation with singularity at $z_0$ is
\[
G_{z_0}(z) = \frac{1}{{(4\pi (t_0  - t))^{m/2} }}\exp \left( {-\frac{{(x - x_0 )^2 }}{{4(t_0  - t)}}} \right)
\]
where $t<t_0$. We also write $G = G_{(0,0)}$. Let $i(M)$ be the injectivity radius of $M$, and suppose that $(\varphi,A)$ is a solution to the flow (\ref{flow1})-(\ref{flow2}) on $M \times [0,T)$. Let $\phi_x$ be a smooth cut-off function with $\left| \phi_x  \right| \leqslant 1$, $\phi_x  \equiv 1$ on $B_{i(M)/2}(x)$, $\phi_x \equiv 0$ outside $B_{i(M)}(x)$ and $\left| \nabla \phi_{x} \right| \leq c/i(M)$ for some constant $c$. We also abbreviate $\phi = \phi_{x_0}$. We define
\begin{equation} \label{capitalphi}
\Phi (R;\varphi ,A) = R^2 \int_{T_R (z_0 )} {e(\varphi ,A)(z)\phi ^2 G}dVdt,
\end{equation}
and
\begin{equation} \label{Ffunction}
\mathscr{F} (R;\varphi ,A)= \int_{T_R(z_0) } {Rt\left( {\left| {\frac{{\partial A}}{{\partial t}} + \frac{{x_k }}{{2t}}\frac{\partial }{{\partial x_k }}\rfloor F_A } \right|^2  + 2\left| {\frac{{\partial \varphi }}{{\partial t}} + \frac{{x_k }}{{2t}}\nabla _A^k \varphi } \right|^2 } \right)\phi ^2 G \sqrt g dz},
\end{equation}
where 
\[
\frac{\partial }{{\partial x_k  }} \rfloor F_A  = F_A (\frac{\partial }
{{\partial x_k  }}, \cdot )=F^{kj}dx^j
\]
defines a 1-form.

\begin{Lemma3} \label{monotonicity}
Let $(\varphi,A)$ be a smooth solution of (\ref{flow1})-(\ref{flow2}) on $M \times [0,T)$ with intial data $(\varphi _0 ,A_0 )$. Then for $z_0 \in M \times [0,T]$ and any $R_a$ and $R_b$ satisfying $0<R_a \leq R_b \leq R_0$ for some $R_0 \leq \min \left\{ {i(M),\sqrt {t_0 } /2} \right\}$, we have
\[
\Phi (R_a ;\varphi ,A) + \int_{R_a }^{R_b } {e^{aR} F(R)} \leqslant e^{c(R_b  - R_a )}\Phi (R_b ;\varphi ,A) + c(R_b^2  - R_a^2 )\mathcal{SW}(\varphi _0,A_0),
\]
where $c$ depends only on the geometry of $M$.
\end{Lemma3}
\emph{Proof}. We show that
\begin{equation} \label{tobeshown}
\frac{d}{{dR}}\Phi (R;\varphi ,A) + \mathscr{F} (R;\varphi ,A)  \geqslant  - c\Phi (R ;\varphi ,A) - cR \mathcal{SW}(\varphi _0 ,A_0 ).
\end{equation}
The required result then follows by multiplying (\ref{tobeshown}) by $e^{aR}$ for some sufficiently large $a>0$, and integrating from $R_a$ to $R_b$. To show (\ref{tobeshown}), we may assume that $z_0 = (0,0)$, which implies that $t<0$ on $T_R$. We rescale the coordinates to $x = R \tilde x$, $t = R^2 \tilde t$. In these coordinates,
\[
\Phi (R;\varphi , A) = \int_{T_1 } {R^4 e(\varphi ,A)(R\tilde x,R^2 \tilde t)\phi ^2 (R\tilde x)G(\tilde z)\sqrt {g(R\tilde x)} d\tilde z} 
\]
where $d\tilde z = d\tilde xd\tilde t$. For some $R \leq R_0$, we compute
\begin{align*}
\frac{d}{{dR}}\Phi (R;\varphi ,A) & = \int_{T_1 } {\frac{d}
{{dR}}\left[ {R^4 e(\varphi ,A)(R\tilde x,R^2 \tilde t)\phi ^2 (R\tilde x)\sqrt {g(R\tilde x)} } \right]G(\tilde z) d\tilde z} \\
& =  \int_{T_1 } {4R^3 e(\varphi ,A)(R\tilde x,R^2 \tilde t)\phi ^2 (R\tilde x) G(\tilde z)\sqrt {g(R\tilde x)} d\tilde z} \\
& + \int_{T_1 } {R^4 \tilde x_k \frac{\partial }
{{\partial x_k }}e(\varphi ,A)(R\tilde x,R^2 \tilde t)\phi ^2 (R\tilde x) G(\tilde z)\sqrt {g(R\tilde x)} d\tilde z}  \\
& + \int_{T_1 } {2R^5\tilde t\frac{\partial }
{{\partial t}}e(\varphi ,A)(R\tilde x,R^2 \tilde t)\phi ^2 (R\tilde x) G(\tilde z)\sqrt {g(R\tilde x)} d\tilde z} \\
& + \int_{T_1 } {R^4 e(\varphi ,A)(R\tilde x,R^2 \tilde t)\tilde x_k \frac{\partial }
{{\partial x_k }}\left( {\phi ^2 \sqrt g } \right)(R\tilde x)G(\tilde z)d\tilde z}  \\
& := I_1  + I_2  + I_3  + I_4. 
\end{align*}
Rescaling coordinates back to (x,t), we have
\[
I_1=  \int_{T_R } {4R e(\varphi ,A)\phi ^2 G\sqrt g dz}
\]
and
\[
I_4  = \int_{T_R } {Re(\varphi ,A)} x_k \frac{\partial }
{{\partial x_k }}\left( {\phi ^2 \sqrt g } \right)Gdz.
\]
For the second term,
\[
I_2  = \int_{T_R } {Rx_k \frac{\partial }{{\partial x_k }}e(\varphi ,A)\phi ^2 G\sqrt g dz}. 
\]
This simplifies as follows:
\[
\frac{\partial }
{{\partial x_k }}\left[ {\left| {\nabla _A \varphi } \right|^2+\frac{1}{2}\left| {F_A } \right|^2+ \left( {\frac{S}
{4}\left| \varphi  \right|^2  + \frac{1}
{8}\left| \varphi  \right|^4 } \right)} \right]
\]
\[
=\left\langle {\nabla _M^k F_A ,F_A } \right\rangle  + 2\operatorname{Re} \left\langle {\nabla _A^k \nabla _A^j \varphi ,\nabla _A^j \varphi } \right\rangle  + \frac{1}
{2}\left( {S + \left| \varphi  \right|^2 } \right)\operatorname{Re} \left\langle {\nabla _A^k \varphi ,\varphi } \right\rangle.
\]
Note that
\[
2\operatorname{Re} \left\langle {\nabla _A^k \nabla _A^j \varphi ,\nabla _A^j \varphi } \right\rangle  = 2\operatorname{Re} \left\langle {\nabla _A \nabla _A^k \varphi ,\nabla _A \varphi } \right\rangle  - 2\operatorname{Re} \left\langle {\Omega _A^{jk} \varphi ,\nabla _A^j \varphi } \right\rangle. 
\]
Using the fact that
\[
\frac{{\partial G}}{{\partial x_j }} = \frac{{x_j }}{{2t}}G,
\]
we have
\[
 - 2\int_{T_R } {Rx_k \operatorname{Re} \left\langle {d( G ) \wedge \nabla _A^k \varphi ,\nabla _A \varphi } \right\rangle } \phi ^2\sqrt g dz =  - 4\int_{T_R } {Rt\left| {\frac{{x_k }}{{2t}}\nabla _A^k \varphi } \right|^2 } \phi ^2 G \sqrt g dz.
\]
For the curvature term, we recall the Bianchi identiy $dF_A=0$ which implies that
\[
\partial _k F^{ij}  = \partial _i F^{kj}  - \partial _j F^{ki}, 
\]
and we compute in local coordinates
\begin{align*}
x_k \partial _k \sum\limits_{i < j} {(F^{ij} )^2 G\phi ^2 } & = 2x_k \sum\limits_{i < j} {F^{ij} \partial _k F^{ij} } G\phi ^2 \\
& = 2x_k \sum\limits_{i<j} {F^{ij} (\partial _i F^{kj}  - \partial _j F^{ki} )} G\phi ^2 \\
& = x_k F^{ij} \partial _i F^{kj} G\phi ^2 \\
& = \partial _i (x_k F^{ij} F^{kj} G\phi ^2 ) - x_k F^{kj} \partial _i F^{ij} G\phi ^2  - F^{kj} F^{ij} \partial _i (x_k G\phi ^2 ).
\end{align*}
Observe that the first term will integrate to zero by Stoke's theorem, and for the second term we have $(d^* F_A )_j  = \partial _i F^{ij}$. To deal with the third term, we see that
\[
-x_k F^{kj} x_i F^{ij} \frac{1}{{2t}}G \phi ^2 = -2t\left| {\frac{{x_k }}{{2t}}\frac{\partial }{{\partial x_k }} \rfloor F_A } \right|^2 G \phi ^2. 
\]
Note that $\left| {d\phi } \right|G \leqslant c$ since $\left| {d\phi } \right| = 0$ on $B_{i(M)/2}(x_0)$. Then
\begin{align*}
I_2 & = \int_{T_R } {Rx_k \left[ {\frac{1}{2} \partial _k \left| {F_A } \right|^2 + 2\operatorname{Re} \left\langle {\nabla _A \nabla _A^k \varphi ,\nabla _A \varphi } \right\rangle  - 2\operatorname{Re} \left\langle {\Omega _A^{jk} \varphi ,\nabla _A^j \varphi } \right\rangle } \right.} \\
&\left. { + \frac{1}{2}\left( {S + \left| \varphi  \right|^2 } \right)\operatorname{Re} \left\langle {\nabla _A^k \varphi ,\varphi } \right\rangle } \right] \phi ^2 G\sqrt g dz\\
&= \int_{T_R } {Rx_k \left[ {\frac{1}{2} \partial _k \left| {F_A } \right|^2 + 2\operatorname{Re} \left\langle {\nabla _A^k \varphi ,\nabla _A^* \nabla _A \varphi } \right\rangle } \right.} \\
& - \left. {2\operatorname{Re} \left\langle {\Omega _A^{jk} \varphi ,\nabla _A^j \varphi } \right\rangle  + \frac{1}
{2}\left( {S + \left| \varphi  \right|^2 } \right)\operatorname{Re} \left\langle {\nabla _A^k \varphi ,\varphi } \right\rangle } \right]\phi ^2 G \sqrt g dz \\
& - 2\int_{T_R } {Rx_k \operatorname{Re} \left\langle {d(\phi ^2 Gx_k ) \nabla _A^k \varphi ,\nabla _A \varphi } \right\rangle } \sqrt g dz \\
& \geqslant  - \int_{T_R } {R\left| x \right|\left( {\left| {\frac{{\partial A}}
{{\partial t}}} \right|\left| {F_A } \right| + 2\left| {\frac{{\partial \varphi }}
{{\partial t}}} \right|\left| {\nabla _A \varphi } \right|} \right)} \phi ^2 G \sqrt g dz \\
& - 4\int_{T_R } {Rt\left| {\frac{{x_k }}{{2t}}\nabla _A^k \varphi } \right|^2 } \phi ^2 G \sqrt g dz - \int_{T_R } R t\left| {\frac{{x_k }}{{2t}}\frac{\partial }{{\partial x_k }}F_A } \right|^2 G\phi ^2 \sqrt g dz\\
&- c\Phi (R ;\varphi ,A) - cR \mathcal{SW}(\varphi _0 ,A_0 ),
\end{align*}
where we note that $\Omega _A  = \Omega _{A_0 }  + \frac{1}{2}F_A$, and we also recall the fact from \cite{Schabrun.Flow} that
\[
\operatorname{Re} \left\langle {\frac{{\partial A}}
{{\partial t}}\varphi ,\nabla _A \varphi } \right\rangle
= \left\langle {\frac{{\partial A}}
{{\partial t}},i\operatorname{Im} \left\langle {\nabla _A \varphi
,\varphi } \right\rangle } \right\rangle;
\]
note that $\frac{\partial A}{\partial t}$ can be replaced by any 1-form. For the third term,
\begin{align*}
I_3 & = \int_{T_R } {2Rt \left[ {\left\langle {d\frac{{\partial A}}
{{\partial t}},F_A } \right\rangle  + 2\operatorname{Re} \left\langle {\nabla _A \frac{{\partial \varphi }}
{{\partial t}},\nabla _A \varphi } \right\rangle  + \operatorname{Re} \left\langle {\frac{{\partial A}}
{{\partial t}}\varphi ,\nabla _A \varphi } \right\rangle } \right.} \\
& \left. { + \left( {\frac{1}
{2}\left[ {S + \left| \varphi  \right|^2 } \right]\operatorname{Re} \left\langle {\frac{{\partial \varphi }}
{{\partial t}},\varphi } \right\rangle } \right)} \right]\phi ^2 G\sqrt {g}dz\\
& =  - \int_{T_R } {2Rt\left[ {\left| {\frac{{\partial A}}
{{\partial t}}} \right|^2  + 2\left| {\frac{{\partial \varphi }}
{{\partial t}}} \right|^2 } \right]\phi ^2 G\sqrt {g} dz} \\
& - \int_{T_R } {2Rt\left\langle {d\left( {\phi ^2 G} \right) \wedge  \frac{{\partial A}}
{{\partial t}},F_A } \right\rangle \sqrt {g} dz} \\
& - \int_{T_R } {4Rt\operatorname{Re} \left\langle { d\left( {\phi ^2G} \right) \wedge \frac{{\partial \varphi }}
{{\partial t}} ,\nabla _A \varphi } \right\rangle\sqrt {g} dz}.
\end{align*}
Next we obtain
\[
\left\langle {dG  \wedge \frac{{\partial A}}
{{\partial t}},F_A } \right\rangle  = \left\langle {\frac{{x_k }}
{{2t}}dx^k  \wedge \frac{{\partial A}}
{{\partial t}},F_A } \right\rangle G  = \left\langle {\frac{{\partial A}}
{{\partial t}},\frac{{x_k }}
{{2t}}\frac{\partial }
{{\partial x_k }}\rfloor F_A } \right\rangle G, 
\]
and
\[
\left\langle {dG  \wedge \frac{{\partial \varphi }}
{{\partial t}},\nabla _A \varphi } \right\rangle  = \left\langle {\frac{{x_k }}
{{2t}}dx^k  \wedge \frac{{\partial \varphi }}
{{\partial t}},\nabla _A \varphi } \right\rangle G  = \left\langle {\frac{{\partial \varphi }}
{{\partial t}},\frac{{x_k }}
{{2t}}\nabla _A^k \varphi } \right\rangle G. 
\]
Thus
\begin{align*}
I_3 & \geqslant  - \int_{T_R } {2Rt\left[ {\left| {\frac{{\partial A}}
{{\partial t}}} \right|^2  + 2\left| {\frac{{\partial \varphi }}
{{\partial t}}} \right|^2 } \right]\phi ^2 G \sqrt g dz} \\
& - \int_{T_R } {2Rt\left\langle {\frac{{\partial A}}{{\partial t}},\frac{{x_k }}{{2t}}\frac{\partial }{{\partial x_k }}\rfloor F_A } \right\rangle \phi ^2 G \sqrt g dz}  + c\int_{T_R } {2Rt\left| {\frac{{\partial A}}{{\partial t}}} \right|\left| {F_A } \right|\phi \sqrt g dz} \\
& - \int_{T_R } {4Rt\operatorname{Re} \left\langle {\frac{{\partial \varphi }}{{\partial t}},\frac{{x_k }}
{{2t}}\nabla _A^k \varphi } \right\rangle \phi ^2 G \sqrt g dz}  + c\int_{T_R } {4Rt\left| {\frac{{\partial \varphi }}
{{\partial t}}} \right|\left| {\nabla _A \varphi } \right|\phi \sqrt g dz}. \\
& \geq - \int_{T_R } {Rt\left( {\left| {\frac{{\partial A}}{{\partial t}} + \frac{{x_k }}{{2t}}\frac{\partial }{{\partial x_k }}\rfloor F_A } \right|^2  + 2\left| {\frac{{\partial \varphi }}{{\partial t}} + \frac{{x_k }}{{2t}}\nabla _A^k \varphi } \right|^2 } \right)\phi ^2 G \sqrt g dz} \\
& + \int_{T_R } {Rt\left( {\left| {\frac{{x_k }}{{2t}}\frac{\partial }
{{\partial x_k }}\rfloor F_A } \right|^2  + 2\left| {\frac{{x_k }}{{2t}}\nabla _A^k \varphi } \right|^2 } \right)\phi ^2 G \sqrt g dz} \\
& - \int_{T_R } {Rt\left[ {\left| {\frac{{\partial A}}
{{\partial t}}} \right|^2  + 2\left| {\frac{{\partial \varphi }}
{{\partial t}}} \right|^2 } \right]\phi ^2 G \sqrt g dz} \\
& - cR \mathcal{SW}(\varphi _0 ,A_0 )
\end{align*}
using Young's inequality, the energy inequality (\ref{energyinequality2}), and noting that $\left|t\right| \leq 4R^2$ and $R \leq R_0$. Finally, since as in \cite{CS}, $R^{ - 1} \left| x \right|^2 G \leqslant c(1 + G)$, combining the working above (and recalling that $t<0$ on $T_R$), one obtains (\ref{tobeshown}). $\square$ \vspace{5mm}

\begin{Corollarymono} \label{corollarymono}
There exists a constant $a>0$ such that
\[
e^{aR}\Phi (R ;\varphi ,A) + cR^2 \mathcal{SW}(\varphi _0 ,A_0 )
\]
is non-decreasing in $R$, where $c$ here represents the same constant which appears in (\ref{tobeshown}).
\end{Corollarymono}
\emph{Proof.} The result f-ollows from (\ref{tobeshown}) by multiplying by $e^{aR}$ for some sufficiently large $a>0$, and integrating from $R_a$ to $R_b$.
 $\square$ \vspace{5mm}

\begin{Lemma4} \label{regularitytheorem}
Suppose $(\varphi, A)\in C^{\infty}( P_R (y,s))$ satisfies
(\ref{flow1})-(\ref{flow2}). Then there exist  constants
$\delta$ and $R_1$ such that if $R \leqslant R_1$ and
\[
\mathop {\sup }\limits_{0 <t < s } R^{4-m} \int_{B_R (y)} {\left(
{\left| {\nabla _A \varphi } \right|^2  + \left| {F_A } \right|^2
} \right)\,dV} < \delta,
\]
then
\[
\mathop {\sup }\limits_{P_{R/2} (y,s)} \left( {\left| {\nabla _A
\varphi } \right|^2  + \left| {F_A } \right|^2 } \right) \leqslant
256R^{ - 4}.
\]
\end{Lemma4}
\emph{Proof.} We begin by choosing $r_0<R$ so that
\begin{equation} \label{r0def}
(R - r_0 )^4 \mathop {\sup }\limits_{P_{r_0 } (y,s)} \left(
{\left| {\nabla _A \varphi } \right|^2  + \left| {F_A } \right|^2
} \right) = \mathop {\max }\limits_{0 \leqslant r \leqslant R}
\left[ {(R - r)^4 \mathop {\sup }\limits_{P_r (y,s)} \left(
{\left| {\nabla _A \varphi } \right|^2  + \left| {F_A } \right|^2
} \right)} \right].
\end{equation}
Let
\[
e_0  = \mathop {\sup }\limits_{P_{r_0 } (y,s)} \left( {\left|
{\nabla _A \varphi } \right|^2  + \left| {F_A } \right|^2 }
\right) = \left( {\left| {\nabla _A \varphi } \right|^2  + \left|
{F_A } \right|^2 } \right)(x_0 ,t_0 )
\]
for some $(x_0 ,t_0 ) \in \bar P_{r_0 } (y,s)$. We claim that
\begin{equation} \label{eclaim}
e_0  \leqslant 16(R - r_0 )^{ - 4}.
\end{equation}
Then
\begin{align*}
&(R - r)^4 \mathop {\sup }\limits_{P_r (y,s)} \left( {\left|
{\nabla _A \varphi } \right|^2  + \left| {F_A } \right|^2 }
\right) \leqslant (R - r_0 )^4 \mathop {\sup }\limits_{P_{r_0 }
(y,s)} \left( {\left| {\nabla _A \varphi } \right|^2  + \left|
{F_A } \right|^2 } \right)
\\
&
 \leqslant 16(R - r_0 )^4 (R - r_0 )^{ - 4}  = 16
\end{align*} for any $r<R$. Choosing $r = \frac{1}{2}R$ in the
above, we have the required result.

We now prove (\ref{eclaim}). Define $\rho _0  = e_0^{ - 1/4}$ and
suppose by contradiction that $\rho _0  \leqslant \frac{1}{2}(R -
r_0 )$. We rescale variables $x=x_0+\rho_0\tilde x$ and
$t=t_0+\rho_0^2 \tilde t$ and set
\[
\psi (\tilde x,\tilde t) = \varphi (x_0  + \rho _0 \tilde x,t_0  +
\rho _0^2 \tilde t),
\]
\[
B(\tilde x,\tilde t) = \rho _0 A(x_0  + \rho _0 \tilde x,t_0  +
\rho _0^2 \tilde t),
\]
giving
\[
\left| {\nabla _B \psi } \right|^2  = \rho _0^2 \left| {\nabla _A
\varphi } \right|^2,
\]
\[
\left| {F_B } \right|^2  = \rho _0^4 \left| {F_A } \right|^2.
\]
We define
\[
e_{\rho _0 } (\tilde x,\tilde t) = \left| {F_B } \right|^2  + \rho
_0^2 \left| {\nabla _B \psi } \right|^2  = \rho _0^4 \left(
{\left| {\nabla _A \varphi } \right|^2  + \left| {F_A } \right|^2
} \right)
\]
so that
\[
e_{\rho _0 } (\tilde x,\tilde t) \leqslant e_{\rho _0 } (0,0) = 1.
\]
We compute
\begin{align*}
\mathop {\sup }\limits_{\tilde P_1 (0,0)} e_{\rho _0 } (\tilde
x,\tilde t)
& = \rho _0^4 \mathop {\sup }\limits_{P_{\rho _0 } (x_0 ,t_0 )} \left( {\left| {\nabla _A \varphi } \right|^2  + \left| {F_A } \right|^2 } \right) \\
& \leqslant \rho _0^4 \mathop {\sup }\limits_{P_{\frac{{R + r_0 }}
{2}} (y,s)} \left( {\left| {\nabla _A \varphi } \right|^2  + \left| {F_A } \right|^2 } \right)\\
& = \rho _0^4 \left( {\frac{{R - r_0 }} {2}} \right)^{ - 4} \left(
{R - \frac{{R + r_0 }} {2}} \right)^4 \mathop {\sup
}\limits_{P_{\frac{{R + r_0 }}
{2}} (y,s)} \left( {\left| {\nabla _A \varphi } \right|^2  + \left| {F_A } \right|^2 } \right)\\
& \leqslant \rho _0^4 \left( {\frac{{R - r_0 }}
{2}} \right)^{ - 4} \left( {R - r_0 } \right)^4 e_0  = 16,
\end{align*}
where we have used that $P_{\rho _0 } (x_0 ,t_0 ) \subset
P_{\frac{{R + r_0 }} {2}} (y,s)$, and to get to the last line we
have used (\ref{r0def}). This implies that
\[
e_{\rho _0 }  = \rho _0^4 \left( {\left| {\nabla _A \varphi }
\right|^2  + \left| {F_A } \right|^2 } \right) \leqslant 16
\]
on $\bar P_1(0,0)$. By Lemma \ref{firstderivativeestimate},
\begin{align*}
(\frac{\partial } {{\partial t}}+\Delta )  \left( {\left| {\nabla
_A \varphi } \right|^2  + \left| {F_A } \right|^2 }  +1 \right)
\leqslant  c\left( {\left| {F_A } \right| + 1} \right) \left( {
{\left| {\nabla _A \varphi } \right|^2  + \left| {F_A } \right|^2
}   + 1} \right).
\end{align*}
Then
\begin{align*}
\left( {\frac{\partial }{{\partial \tilde t}} + \tilde \Delta }
\right) (e_{_{\rho _0 } }+\rho_0^4)
& = \rho _0^6 \left( {\frac{\partial }{{\partial t}} + \Delta } \right)\left( {\left| {\nabla _A \varphi } \right|^2  + \left| {F_A } \right|^2 } \right)\\
& \leqslant c\rho _0^6 \left( {\left| {F_A } \right| + 1}
\right)\left( {\left| {\nabla _A \varphi } \right|^2  + \left|
{F_A } \right|^2 }+1 \right)
\end{align*}
on $\bar P_1(0,0)$. Note that by assumption $\rho _0  < R$,
$\rho_0^2 |F_A|$ is thus bounded by a constant. Then
\[
\left( {\frac{\partial } {{\partial \tilde t}} + \tilde \Delta }
\right)\left( {e_{\rho _0 }  +  \rho_0^4} \right) \leqslant
c\left( {e_{\rho _0 }  + \rho_0^4} \right)
\]
for a constant $c>0$.  We apply Moser's Harnack inequality to give

\begin{align*}
1 + \rho_0^4 = e_{\rho _0 } (0,0) + \rho_0^4
& \leqslant c\int_{\tilde P_1 (0,0)} {e_{\rho _0 } d\tilde xd\tilde t} +c\rho_0^4 \\
& = c\rho _0^{2-m} \int_{P_{\rho _0 } (x_0 ,t_0 )} {\left( {\left| {\nabla _A \varphi } \right|^2  + \left| {F_A } \right|^2 } \right)} dVdt  +c\rho_0^4\\
& \leqslant c\mathop {\sup }\limits_{0  \leqslant t \leqslant s } R^{4-m}
\int_{B_R (y)} {\left( {\left| {\nabla _A \varphi } \right|^2  +
\left| {F_A } \right|^2 } \right)}dV
 +  cR^4 \\
& < c\delta +cR^4,
\end{align*}
where  we have used that $\rho _0  < R$. Now if we choose $R_1$
 and $\delta $  sufficiently small, we have
the desired contradiction. $\square$ \vspace{5mm}

\section{Singularity analysis}

Let $(\varphi,A)$ be a smooth solution on $[0,T)$. Suppose there exists some $R \leq R_1$ such that
\[
R^{4-m} \int_{B_R (y)} {\left({\left| {\nabla _A \varphi } \right|^2  + \left| {F_A } \right|^2
} \right)\,dV} < \delta,
\]
for all $x_0 \in M$ and $t_0 = T$. Then by Lemma \ref{regularitytheorem}, $\left| {\nabla _A \varphi } \right|^2$ and $\left| {F_A } \right|^2$ are uniformly bounded on $M \times [0,T)$. As in \cite{Schabrun.Flow}, Using Lemma \ref{higherderivativeestimate} we can show that $\phi$ and $A$ are smooth at $t=T$. In conjunction with the local existence result, we can extend $(\varphi,A)$ to a global smooth solution.

We define the singular set
\[
\sum  = \bigcap\limits_{0 < R \leq R_1} {\left\{ {x_0  \in M: \mathop {\lim \sup }\limits_{t \to T} R^{4-m} \int_{B_{_R } (x_0)} {\left( {\left| {\nabla _A \varphi } \right|^2  + \left| {F_A } \right|^2 } \right)dV \geq \delta } } \right\}}.
\]
By the above discussion, $(\varphi(T),A(T))$ is smooth on $M\backslash \Sigma$. Let $\Sigma '$ be defined as for $\Sigma$, but with $\varepsilon_0$ replaced by a smaller constant. Clearly $\Sigma  \subseteq \Sigma '$. Furthermore, if $x \in M \backslash \Sigma$ then by smoothness at $x$, $x \in M \backslash \Sigma '$. Thus replacing $\varepsilon_0$ with with any smaller constant defines the same set.

If $x \in M \backslash \Sigma$, then by Lemma \ref{regularitytheorem}, $B_R(x) \in M \backslash \Sigma$ for some $R$. Thus $\Sigma$ is closed. Unlike in the 4 dimensional case \cite{Schabrun.Flow}, we cannot conclude at this point that the singular set is finite. Although it will not be needed in the proof of Theorem 1, one can instead show that $\Sigma$ has finite $(m-4)$-dimensional Hausdorff measure $\mathcal{H}^{m-4}$. Explicitly, for $x_0 \in \Sigma$,
\begin{equation} \label{Rmestimate}
\delta  \leq  \mathop {\lim \sup }\limits_{t \to T} R^{4 - m} \int_{B_R (x_0 )} {e(\varphi ,A)}dV
\end{equation}
for any $R$. The family $\mathcal{F} = \{ B_R (x_0 ):x_0  \in \Sigma \}$ covers $\Sigma$, and by Vitali's covering lemma, there exists a finite subfamily $\mathcal{F'} = \{ B_R (x_j )\}$ such that any two balls in $\mathcal{F'}$ are disjoint and $\{B_{5R} (x_j )\}$ covers $\Sigma$. Then using (\ref{Rmestimate}),
\begin{align*}
\sum\limits_j {(5R)^{m-4} } & \leqslant \frac{{5^m }}{{\delta }}\sum\limits_j {\mathop {\lim \sup }\limits_{t \to T} \int_{B_R(x_j)} {\left( {\left| {\nabla _A \varphi } \right|^2  + \left| {F_A } \right|^2 } \right)dV} } \\
&  \leqslant C \mathcal{SW}(\varphi _0 ,A_0 ),
\end{align*}
where $\{B_{5R} (x_j )\}$ covers $\Sigma$. It follows that $\mathcal{H}^{m-4}(\Sigma)$ is finite, as claimed.

To establish Theorem 1, we show that $\Sigma  = \emptyset$. Suppose by contradiction that $\Sigma$ is non-empty. Since the flow is smooth on $[0,T)$, we can find sequences $x_n \in M$, $t_n \to T$, $R_n \to 0$ such that
\begin{align} \label{blowupboundhigher}
\delta & > R_n^{4-m} \mathcal{SW}_{B_{R_n } (x_n )} (\varphi (t_n ),A(t_n
))\nonumber\\
& = \mathop {\sup }\limits_{0  \leqslant t \leqslant t_n ,\,\,x
\in M}R_n^{4-m} \mathcal{SW}_{B_{R_n } (x)} (\varphi
(t),A(t))>\frac {\delta} {2}
\end{align}
for each $n$, where $\mathcal{SW}_{B_R(x)}$ is defined by
\[
\mathcal{SW}_{B_R (x)} (\varphi ,A) = \int_{B_R (x)} {\left|
{\nabla _A \varphi } \right|^2  + \frac{1} {2}\left| {F_A }
\right|^2  + \frac{S} {4}\left| \varphi  \right|^2  + \frac{1}
{8}\left| \varphi  \right|^4 }.
\]
By the compactness of $M$, passing to a subsequence we have $x_n \to x_0$ where $x_0  \in \Sigma$ by Lemma \ref{regularitytheorem}. We define the region
\[
D_n  = \left\{ {(y,s):R_n y + x_n  \in B_{i(M)/2}(x_n),s \in \left[ { -R_n^{-2}t_n, 0} \right]} \right\}= : U_n \times \left[ { -R_n^{-2}t_n, 0} \right].
\]
Note that as $n \to \infty$, $D_n  \to \mathbb{R}^m  \times \left( { -\infty, 0} \right]$. Furthermore, truncating the sequence if necessary, we can arrange that $B_{i(M)/2}(x_n) \subset B_{i(M)}(x_0)$. We rescale to
\[
\varphi _n (y,s) = \varphi (R_n y + x_n ,R_n^2 s + t_n ),
\]
\[
A_n (y,s) = R_n A(R_n y + x_n ,R_n^2 s + t_n ).
\]
which are defined on $D_n$. We have
\[
\left| {\nabla _{A_n } \varphi _n } \right|^2  = R_n^2 \left| {\nabla _A \varphi } \right|^2, 
\]
\[
\left| {F_{A_n } } \right|^2  = R_n^4 \left| {F_A } \right|^2. 
\]
If we choose our local coordinates on $B_{i(M)}(x_0)$ to be orthonormal coordinates, then the metric on the rescaled space is simply $g_{ij}=\delta_{ij}$. From (\ref{blowupboundhigher}),
\begin{equation} \label{tocontradict}
 \int_{B_1(0)} {{R_n^2 \left| {\nabla _{A_n } \varphi _n } \right|^2  + \left| {F_{A_n } } \right|^2  + R_n^4 \left( {\frac{S}{4}\left| \varphi_n  \right|^2  + \frac{1}{8}\left| \varphi_n  \right|^4 } \right)}dy > \frac{{\delta }}{2}}
\end{equation}
for each $n$ and $s=0$. Next, from Lemma \ref{regularitytheorem} and (\ref{blowupboundhigher}) we have
\begin{equation} \label{scaledbounds}
\mathop {\sup }\limits_{D_n} \left( {\left| {\nabla _{A_n } R_n \varphi _n } \right|^2  + \left| {F_{A_n } } \right|^2 } \right) \leqslant K_1, 
\end{equation}
where $K_1$ is independent of $n$. We consider the rescaled equations
\begin{equation} \label{phirescaledequation}
\frac{{\partial R_m \varphi _m }} {{\partial s}} = R_m^3
\frac{{\partial \varphi }} {{\partial t}} =  - \nabla _{A_m }^*
\nabla _{A_m } R_m \varphi _m  - \frac{1} {4}\left[ {R_m^2 S +
\left| {R_m \varphi _m } \right|^2 } \right]R_m \varphi _m,
\end{equation}
\begin{equation} \label{Arescaledequation}
\frac{{\partial A_m }} {{\partial s}} = R_m^3 \frac{{\partial A}}
{{\partial t}} =  - d^* F_{A_m }  - i\operatorname{Im}
\left\langle {\nabla _{A_m } R_m \varphi _m ,R_m \varphi _m }
\right\rangle.
\end{equation}
Noting the similarity of these equations to
(\ref{flow1}) and (\ref{flow2}), we use (\ref{scaledbounds})
and results identical to Lemma (\ref{higherderivativeestimate}) to find
\begin{equation} \label{higherblowupbounds}
\mathop {\sup }\limits_{D_n} \left(\left| {\nabla _{A_n }^{(k+1)} R_n \varphi _n } \right|^2  +
\left| {\nabla _M^{(k)} F_{A_n } } \right|^2 \right) \leqslant K_{k+1}
\end{equation}
for each $k \geq 0$. Thus by a result of Uhlenbeck (theorem 1.3 in \cite {U2}, see also \cite{Hong.Tian}), passing to a subsequence and using an appropriate gauge we have $C^\infty$ convergence $R_n \varphi _n  \to \tilde \varphi = 0$ (since $\varphi_n$ is bounded), $A_n  \to \tilde A$ where $\tilde \varphi$ and $\tilde A$ are defined on $\mathbb{R}^m  \times \left( { -\infty, 0} \right]$. Then as $n \to \infty$ in (\ref{tocontradict}), we obtain
\begin{equation} \label{tocontradict2}
\int_{B_2(0)} {\left| {F_{\tilde A} } \right|^2 dy \geq \frac{\delta}{2} }
\end{equation}
for $s=0$. Since $R_n \varphi_n \to 0$, from (\ref{Arescaledequation}) $\tilde A$ satisfies the equation
\[
\frac{{\partial \tilde A }} {{\partial s}} = - d^* F_{\tilde A}
\]
on $\mathbb{R}^m  \times \left( { -\infty, 0} \right]$. Using the Bianchi identiy $dF_{\tilde A}=0$, This implies that
\[
\frac{{\partial F_{\tilde A} }} {{\partial s}} = - \Delta F_{\tilde A}
\]
on $\mathbb{R}^m  \times \left( { -\infty, 0} \right]$, where $\Delta  = d^* d + dd^*$. Since the solution to the heat equation converges to constant data in infinite time, the only possible solution to this equation satisfying (\ref{scaledbounds}) is $F_{\tilde A} = $ constant. See for example theorem 9 of chapter 2 in \cite{Evans}. In the notation of  \cite{Evans}, choose $k=1$ and $t=0$, note that for us $\left\| u \right\|_{L^1 (C(x,t;r))}  \leqslant cr^{n + 2}$, and let $r \to \infty$.

\section{Proof of Theorem 1}

As in \cite{TianCalibrated} and \cite{LinWang}, the term $\mathscr{F} (R;\varphi ,A)$ in Lemma \ref{monotonicity} can be used to further analyse the singularity (see for example Lemma 3.3.2 of \cite{TianCalibrated}). However, we are already in a position to show that the existence of a singularity implies a contradiction. Noting that $G_{(x_n ,t_n )}  \geqslant cR_n^{-m}$ on $B_{rR_n } (x_n ) \times \left[ {t_n  - 4(rR_n )^2 ,t_n  - (rR_n )^2 } \right]$, we consider for any $r \in (0,\infty)$,
\begin{align*}
& r^{2 - m} \int_{B_r (0) \times \left[ { - 4r^2 , - r^2 } \right]} {\left| {F_{\tilde A} } \right|^2 } dyds \\ & =\mathop {\lim }\limits_{n \to \infty } (rR_n )^{2 - m} \int_{B_{rR_n } (x_n ) \times \left[ {t_n  - 4(rR_n )^2 ,t_n  - (rR_n )^2 } \right]} {\left| {F_A } \right|^2 } dVdt \\
& \leqslant c \mathop {\lim }\limits_{n \to \infty } (rR_n )^2 \int_{B_{rR_n } (x_n ) \times \left[ {t_n  - 4(rR_n )^2 ,t_n  - (rR_n )^2 } \right]} {\left| {F_A } \right|^2 } G_{(x_n ,t_n )} dVdt \\
& \leqslant c \mathop {\lim }\limits_{n \to \infty } (rR_n )^2 \int_{T_{rR_n } (x_0 ,T)} {e(\varphi ,A)} G_{(x_0 ,T)} \phi ^2 dVdt.
\end{align*}
However, the latter expression is bounded by Lemma \ref{monotonicity}. Thus
\[
\int_{B_r (0) \times \left[ { - 4r^2 , - r^2 } \right]} {\left| {F_{\tilde A} } \right|^2 } dyds \leq cr^{m-2}.
\]
But since $\left| {F_{\tilde A} } \right|$ is constant and non-zero by (\ref{tocontradict2}), this implies that $r^4 \leq c$. this is impossible for $r$ sufficiently large. This proves Theorem 1. $\square$ \vspace{5mm}

\bibliographystyle{plain}

\end{document}